\input amstex
\documentstyle{amsppt}
\pagewidth{12.4cm}
\pageheight{18.4cm}

\def\>{\rightarrow}
\def\<{\leftarrow}

\def\[{\lbrack}
\def\]{\rbrack}

\def\al{\alpha}
\def\be{\beta}

\def\ph{\varphi}

\def\ps{\psi}

\def\Ph{\Phi}

\def\E{\text{End}}

\topmatter
\title Knit Products of Graded Lie Algebras and Groups 
\endtitle
\author Peter W. Michor \endauthor
\affil Institut f\"ur Mathematik \\ Universit\"at Wien \\ Austria \endaffil
\address{Institut f\"ur Mathematik der Universit\"at Wien,
Strudlhofgasse 4, A-1090 Wien, Austria}\endaddress
\abstract{If a graded Lie algebra is the direct sum of two
graded sub Lie algebras, its bracket can be written in a form
that mimics a "double sided semidirect product". It is called
the {\it knit product} of the two subalgebras then. The
integrated version of this is called a {\it knit product} of
groups --- it coincides with the {\it Zappa-Sz\'ep product}. The
behavior of homomorphisms with respect to knit products is investigated.}
\endabstract
\thanks
This paper is in final form and no version of it will
appear elsewhere.
\endthanks

\endtopmatter
\document
\heading Introduction \endheading

If a Lie algebra is the direct sum of two sub Lie algebras one
can write the bracket in a way that mimics semidirect products
on both sides. The two representations do not take values in the
respective spaces of derivations; they satisfy equations (see
1.1) which look "derivatively knitted" --- so we call them a
derivatively knitted pair of representations. These equations
are familiar for the Fr\"olicher-Nijenhuis bracket of
differential geometry, see \cite{1} or \cite{2, 1.10}. This paper is the
outcome of my investigation of what formulas 1.1 mean
algebraically. It was a surprise for me that they describe the
general situation (Theorem 1.3). Also the behavior of
homomorphisms with respect to knit products is investigated
(Theorem 1.4).

The integrated version of a knit product of Lie algebras will be
called a knit product of groups --- but it is well known to
algebraists under the name {\it Zappa-Sz\'ep product}, see \cite{3}
and the references therein. I present
it here with different notation in order to describe afterwards
again the behavior of homomorphisms with respect to this product.
This gives a kind of generalization of the method of induced
representations.

\heading 1. Knit products of graded Lie algebras \endheading

\subheading{1.1. Definition} Let $A$ and $B$ be graded Lie algebras,
whose grading is in $\bold Z$ or $\bold Z_2$, but only one of
them. A {\it derivatively knitted pair of representations}
$(\al,\be)$ for $(A,B)$ are graded Lie algebra homomorphisms
$\al : A \> \E (B)$ and $\be : B \> \E (A)$ such that:
$$\multline \al(a)[b_1,b_2] = [\al(a)b_1,b_2] +
(-1)^{|a||b_1|}[b_1, \al(a)b_2] -\\ 
- \Bigl ((-1)^{|a||b_1|}\al(\be(b_1)a)b_2
-(-1)^{(|a|+|b_1|)|b_2|}\al(\be(b_2)a)b_1 \Bigr )\endmultline $$
$$\multline \be (b)[a_1,a_2]= [\be(b)a_1,a_2] +
(-1)^{|b||a_1|}[a_1,\be(b)a_2] -\\ 
- \Bigl ((-1)^{|b||a_1|}\be(\al(a_1)b)a_2 - 
(-1)^{(|b| + |a_1|)|a_2|}\be(\al(a_2)b)a_1 \Bigr )\endmultline$$
Here $|a|$ is the degree of $a$. For (non-graded) Lie algebras
just assume that all degrees are zero.

\proclaim{1.2. Theorem} Let $(\al,\be)$ be a derivatively knitted pair
of representations for graded Lie algebras $A=\bigoplus A_k$ and
$B=\bigoplus B_k$. Then $A\oplus B := \bigoplus_{k,l}
(A_k\oplus B_l)$ becomes a graded Lie algebra
$A\oplus_{(\al,\be)}B$ with the following bracket:
$$\multline 
[(a_1,b_1),(a_2,b_2)] := \Bigl ([a_1,a_2] + \be(b_1)a_2 -
(-1)^{|b_2||a_1|} \be(b_2)a_1,\\
[b_1,b_2] + \al(a_1)b_2 - (-1)^{|a_2||b_1|}
\al(a_2)b_1 \Bigr ) \endmultline $$
The grading is $(A\oplus B)_k := A_k \oplus B_k$.
\endproclaim

\demo{Proof} Obviously this bracket is graded anticommutative.
The graded Jacobi identity is checked by computation.
\qed\enddemo

We call $A \oplus_{(\al,\be)} B$ the {\it knit product} of $A$
and $B$. If $\be = 0$ then $\al$ has values in the space of
(graded) derivations of $A$ and $A \oplus 0$ is an ideal in $A
\oplus_{(\al,0)} B$ and we get a semidirect product of graded Lie
algebras. Note also that 
$[(a,0),(0,b)] = \bigl ((-1)^{|b||a|}\be(b)a, \al(a)b
\bigr )$. This is the key to the following theorem.

\proclaim{1.3. Theorem} Let $A$ and $B$ be graded Lie
subalgebras of a graded Lie algebra $C$ such that $A+B=C$ and
$A\cap B=0$. Then $C$ as graded Lie algebra is isomorphic to a
knit product of $A$ and $B$.
\endproclaim

\demo{Proof} For $a\in A$ and $b \in B$ we write 
$$[a,b]=:\al(a)b - (-1)^{|a||b|}\be(b)a$$ 
for the decomposition of $[a,b]$ into components in $C=B+A$. Then
$\be:B\>\E(A)$ and $\al:A\>\E(B)$ are linear. Now decompose both
sides of the graded Jacobi identity
$$[a,[b_1,b_2]]=[[a,b_1],b_2] + (-1)^{|a||b_1|}[b_1,[a,b_2]]$$
and compare the $A$- and $B$-components respectively. This gives
equation 1.1 for $\al$ and that $\be$ is a graded Lie algebra homomorphism. The rest
follows by interchanging $A$ and $B$. Now we decompose 
$[a_1+b_1, a_2+b_2]$ and see that $C=A\oplus_{(\al,\be)}B$.
\qed\enddemo

\subheading{1.4} Now let $\Ph:A\oplus_{(\al,\be)}B \>
A'\oplus_{(\al',\be')}B'$ be a linear mapping between
knit products. Then $\Ph$ can be decomposed into
$\Ph(a,b)=:(f(a)+\ps(b), g(b)+\ph(a))$ for linear mappings
$\ph:A\>B'$, $\ps:B\>A'$, $f:A\>A'$, and $g:B\>B'$.

\proclaim{Theorem} In this situation $\Ph$ is a graded Lie
algebra homomorphism if and only if the following conditions hold:
$$\align 
\ph([a_1,a_2]) &= [\ph(a_1),\ph(a_2)] + \al'(f(a_1))\ph(a_2) \\
&\qquad - (-1)^{|a_1||a_2|}\al'(f(a_2))\ph(a_1) \\
\ps([b_1,b_2]) &= [\ps(b_1),\ps(b_2)] + \be'(g(b_1))\ps(b_2) \\
&\qquad - (-1)^{(|b_1||b_2|} \be'(g(b_2))\ps(b_1) \\
[\ps(b),f(a)] &=  f(\be(b)a) - \be'(g(b))f(a) \\
&\qquad - (-1)^{|a||b|}\bigl (\ps(\al(a)b) - \be'(\ph(a))\ps(b)\bigr) \\
[g(b),\ph(a)] &= \ph(\be(b)a) - \al'(\ps(b))\ph(a) \\
&\qquad -(-1)^{|a||b|} \bigl ( g(\al(a)b) - \al'(f(a))g(b) \bigr )\\
f([a_1,a_2]) &= [f(a_1),f(a_2)] + \be'(\ph(a_1))f(a_2) \\
&\qquad - (-1)^{|a_1||a_2|} \be'(\ph(a_2))f(a_1) \\
g([b_1,b_2]) &= [g(b_1),g(b_2)] + \al'(\ps(b_1))g(b_2) \\
&\qquad -(-1)^{|b_1||b_2|} \al'(\ps(b_2))g(b_1) 
\endalign $$
If $f$ and $g$ are graded Lie algebra homomorphism the last pair of
equations obviously simplifies.
\endproclaim

\demo{Proof} A long but straightforward computation. \qed\enddemo

This theorem can be used to build representations of $C$ out of
representations of $A$ and $B$.

\heading 2. Knit products of groups \endheading

\subheading{2.1. Definition} Let $A$ and $B$ be groups. An {\it
automorphically knitted pair of actions} $(\al,\be)$ for $(A,B)$
are mappings $\al  : B\times A \> A$ and $\be : B\times A \> B$
such that:
\roster
\item $\check \al : B \> \{\text{bijections of A}\}$ is a group
homomorphism, so $\al_{b_1}\circ\al_{b_2} = \al_{b_1b_2}$ and $\al_e =
Id_A$ , where $\al_b(a) := \al(b,a)$.
\item $\check\be : A \> \{\text{bijections of B}\}$ is a group anti
homomorphism, i.e., $\be^{a_1}\circ\be^{a_2} = \be^{a_2a_1}$ and
$\be^e= Id_B$, where $\be^a(b) = \be(b,a)$.
\item $\al_b(a_1a_2) = \al_b(a_1).\al_{\be^{a_1}(b)}(a_2)$.
\item $\be^a(b_1b_2) = \be^{\al_{b_2}(a)}(b_1).\be^a(b_2)$.
\endroster

\proclaim{2.2. Theorem} Let $(\al,\be)$ be an automorphically
knitted pair of actions for $(A,B)$. Then $A\times B$ is a group
$A \times_{(\al,\be)} B$ with the following operations: \newline
\centerline {$(a_1,b_1).(a_2,b_2) := (a_1.\al_{b_1}(a_2),\be^{a_2}(b_1).b_2)$}
\centerline {$(a,b)^{-1} := (\al_{b^{-1}}(a^{-1}),\be^{a^{-1}}(b^{-1}))$.}
Unit is $(e,e)$. $A\times\{e\}$ and $\{e\}\times B$ are
subgroups of $A\times_{(\al,\be)} B$ which are isomorphic to
$A$ and $B$, respectively. If $\check\al \equiv Id_A$
then $\{e\}\times B$ is a normal subgroup of $A
\times_{(\al,\be)}B$ and we have a semidirect product; similarly
if $\check\be \equiv Id_B$.\par
If $A$ and $B$ are topological groups or Lie groups and $\al$,
$\be$ are continuous or smooth, then $A\times_{(\al,\be)}B$ is
also a topological group or Lie group, respectively.
\endproclaim
The proof is routine.

We will call $A\times_{(\al,\be)}B$ the {\it knit product} of
$A$ and $B$ in analogy with section 1. In algebra, with 
different notation, this product is well known under the name
{\it Zappa-Sz\'ep product}. I owe this remark to G. Kowol.

\proclaim{2.3. Theorem} Let $G$ be a group, let $A$ and $B$ be
subgroups such that $G=A.B$ and $A\cap B=\{e\}$. Then $G$ is
isomorphic to a knit product of $A$ and $B$.
\endproclaim

\demo{Proof} Let $b.a = \al(b,a).\be(b,a)$ be the unique
decomposition of $b.a$ in $G=A.B$. Then 
$$a_1b_1a_2b_2 = a_1\al(b_1,a_2)\be(b_1,a_2)b_2 =
(a_1\al_{b_1}(a_2)).(\be^{a_2}(b_1)b_2).$$ 
So it remains to show that $(\al,\be)$ satisfies the conditions
of 2.1. Obviously we have $\al(e,a) = a$, $\be(e,a) = e$,
$\al(b,e) = e$, $\be(b,e) = b$. Comparing coefficients in the
law of associativity of $G$ gives two equations. Setting
suitable elements in these equations to $e$ gives all conditions
of 2.1. 
\qed\enddemo

\subheading{2.4} Let $\Ph = (\Ph_1,\Ph_2) : A\times_{(\al,\be)}B \>
A'\times_{(\al',\be')}B'$ be a mapping between knit products of
groups. We put
$$ \alignat 2 f(a) &:= \Ph_1(a,e),\qquad g(b) &:= \Ph_2(e,b) \tag 1\\
        \ph(b) &:= \Ph_1(e,b),\qquad \ps(a) &:= \Ph_2(a,e)
              \tag 2 \endalignat$$  
Then we have $f:A\>A'$, $g:B\>B'$, $\ph:B\>A'$, $\ps :A\>B'$. 
$\Ph$ is a group homomorphism if and only if 
$$\left\{
\aligned &\Ph_1(a_1\al_{b_1}(a_2), \be^{a_2}(b_1)b_2) = 
           \Ph_1(a_1,b_1).{\al'}_{\Ph_2(a_1,b_1)}(\Ph_1(a_2,b_2))\\
         &\Ph_2(a_1\al_{b_1}(a_2), \be^{a_2}(b_1)b_2) =
           {\be'}^{\Ph_1(a_2,b_2)}(\Ph_2(a_1,b_1)).\Ph_2(a_2,b_2).
     \endaligned\right.\tag3 $$
Now we set in \thetag{3} suitable elements to $e$, use
\thetag{1} and \thetag{2} and
get in turn 
$$\left\{\aligned &\Ph_1(a_1,b_2)= f(a_1).\al'_{\ps(a_1)}(\ph(b_2))\\
                  &\Ph_2(a_1,b_2)= {\be'}^{\ph(b_2)}(\ps(a_1)).g(b_2)
            \endaligned \right. \tag e $$
$$\left\{\aligned &\ph(b_1b_2)= \ph(b_1).\al'_{g(b_1)}(\ph(b_2))\\
                  &\ps(a_1a_2)= {\be'}^{f(a_2)}(\ps(a_1)).\ps(a_2)
            \endaligned \right. \tag f $$
$$\left\{\aligned &\Ph_1(\al_{b_1}(a_2),\be^{a_2}(b_1)) =
                \ph(b_1).\al'_{g(b_1)}(f(a_2)) \\
          &\Ph_2(\al_{b_1}(a_2),\be^{a_2}(b_1)) =
                {\be'}^{f(a_2)}(g(b_1)).\ps(a_2)
            \endaligned \right. \tag 4 $$   
$$\left\{\aligned &f(a_1a_2)= f(a_1).\al'_{\ps(a_1)}(f(a_2))\\
                  &g(b_1b_2)= {\be'}^{\ph(b_2)}(g(b_1)).g(b_2)
            \endaligned \right. \tag g $$
If $f$ and $g$ are homomorphisms of groups then \thetag{g} implies:
$$\left\{\aligned &f(a_2)=  \al'_{\ps(a_1)}(f(a_2))\\
                  &g(b_1)= {\be'}^{\ph(b_2)}(g(b_1))
            \endaligned \right. \tag g' $$
Now we decompose the left hand sides of \thetag{4} with the help of
\thetag{e} and get:
$$\left\{\aligned &f(\al_{b_1}(a_2)).\al'_{\ps(\al_{b_1}(a_2))}
           (\ph(\be^{a_2}(b_1)))= \ph(b_1).\al'_{g(b_1)}(f(a_2)) \\
          &{\be'}^{\ph(\be^{a_2}(b_1))}(\ps(\al_{b_1}(a_2))).
           g(\be^{a_2}(b_1))) = {\be'}^{f(a_2)}(g(b_1)).\ps(a_2)
            \endaligned \right. \tag h $$   

\proclaim{2.5. Theorem} Let $A\times_{(\al,\be)}B$ and
$A'\times_{(\al',\be')}B'$ be knit products of groups and let 
$f:A\>A'$, $g:B\>B'$, $\ph:B\>A'$, $\ps :A\>B'$ be mappings such
that \thetag{f}, \thetag{g}, and \thetag{h} from 2.4 hold. We define 
$\Ph = (\Ph_1,\Ph_2) : A\times_{(\al,\be)}B \>
A'\times_{(\al',\be')}B'$ by 2.4.\thetag{e}, then $\Ph$ is a
homomorphism of groups. 
If $f$ and $g$ are homomorphisms, then  we may use
\thetag{g'} instead of \thetag{g}.
\endproclaim
\demo{Proof} It suffices to check \thetag{3} of 2.5. This is a
difficult computation using 2.4 \thetag{a}-\thetag{h}.\qed\enddemo

For topological groups and Lie groups all the expected
assertions about continuity and smoothness are true. 

This theorem may be used to construct representations of 
$A\times_{(\al,\be)}B$ out of representations of $A$ and $B$ --- a
sort of generalized induced representation procedure.

Starting from the equations 2.1 for a knit product of Lie groups
and deriving the equations of 1.1 for their Lie algebras is a 
very interesting exercise in calculus on Lie groups.

\enddocument